\documentclass[12pt,a4paper]{amsart}
\usepackage{amssymb}
\usepackage{hyperref}

\usepackage{graphicx,amssymb,amsfonts,epsfig,amsthm,a4,amsmath,url}
\usepackage[latin1]{inputenc}

\usepackage{amsmath, amsthm, amssymb, verbatim}
\usepackage{graphicx,amssymb,amsfonts,epsfig,amsthm,a4,amsmath,url,graphicx,enumerate}
\usepackage[latin1]{inputenc}
\usepackage{tikz}
\usepackage{textcomp}
\usetikzlibrary{matrix}
\usepackage{svg}

\newtheorem{thm}{Theorem}[section]
\newtheorem{cor}[thm]{Corollary}

\theoremstyle{definition}

\theoremstyle{remark}
\newtheorem{rem}[thm]{Remark}

\numberwithin{equation}{section}

\begin{document}
\title[Computation of twin-width of graphs]{Computation of twin-width of graphs}
\author[Das]{Kajal Das}
\address{Indian Statistical Institute\\  203 Barrackpore Trunk Road\\Kolkata 700 108}
\email{kdas.math@gmail.com}

\maketitle
\textbf{Abstract:}Twin-width is a recently introduced graph parameter. In this article, we compute twin-width of various finite graphs. In particular, we prove that the twin-widths of finite graphs with 4 and 5 vertices are less than equal to 1 and 2, respectively. We show that the constructions of dual graph and line graph do not preserve twin-width. Also, we give upper bounds for the twin-width of King's graph and Rook's graph.

\textbf{Mathematics Subject Classification (2020):}05C30, 05C38, 05C76, 68R10.

\textbf{Key terms:}:  Graph, Twin-width, Complement graph, Dual graph, Line graph, Graphs with 4 and 5 vertices, King's graph, Rook's graph.

\section{Introduction}

Twin-width is an invariant of graphs introduced in \cite{BKTW20}. It is used to study the parameterized complexity of graph algorithms. It has applications in logic, enumerative combinatorics etc. Recently, it has appeared in many articles (\cite{BGKTW21}, \cite{BGKTW21'}, \cite{BGMSTT21}, \cite{BKRT22}, \cite{BGTT22}, \cite{BCKKLT22}). Moreover, it has been 
studied in the context of finitely generated groups \cite{BGTT22}.  The twin-width is defined for a finite simple graph, later it is extended for a simple infinite graph. The computation of twin-width of a finite graph is extremely difficult. It has been computed before for complete graphs, path graphs, cyclic graphs (or graphs with at most one cycle), Paley graphs, Caterpillar tree, planar graphs etc.  In this article, we compute twin-width for all graphs with vertices 4 and 5 and prove the following theorems. 

\begin{thm}\label{4vertices}
The twin-width of a graph with 4 vertices is less than equal to 1.
\end{thm}

\begin{thm}\label{5vertices}
The twin-width of a graph with 5 vertices is less than equal to 2.
\end{thm}
It is an open problem to determine whether there is an $n$-vertex graph having twin-width at least
$n/2$( see \cite{AHKO22}, page 3) . Therefore, the above-mentioned theorems show that we should look for $n\geq 3$. It is known that twin-with is invariant under taking complement graphs. It was not known how twin-width behaves under other graph operations. In this article, we prove that they are not invariant 
under taking dual graphs and line graphs.

\begin{thm}\label{dual}
Let $\mathcal{C}$ be the collection of simple connected planar graphs whose dual are also a simple connected planar graphs. 
The construction of dual graphs does not preserve twin-width, i.e., there exists a graph $G$ in $\mathcal{C}$ such that 
the twin-width of $G$ and $G^*$ are not equal, where $G^*$ is the dual graph of $G$. 
\end{thm}

\begin{thm}\label{linegraph}
The construction of line graph does not preserve  twin-width.
\end{thm}

However, there are some graphs, King's graph and Rook's graph,  which are associated with Chess 
and are important in Graph Theory. In this article, we study twin-width of these graphs. We briefly recall the definitions of King's graph and Rook's graph. A \textit{King's graph} is a graph that represents all legal moves of the king chess piece on a chessboard where each vertex represents a square on a chessboard and each edge is a legal move. More specifically, an $(n\times m)$-King's graph is a King's graph of an $(n\times m)$-chessboard. On the other hand, a \textit{Rook's graph} is a graph that represents all legal moves of the rook chess piece on a chessboard. Each vertex of a rook's graph represents a square on a chessboard, and each edge connects two squares on the same row or on the same column  (each edge connects the squares that a rook can move between). In this article, we prove the following two theorems.

\begin{thm}\label{King'sgraph}
The twin-width of a $(n\times m)$-King's graph is less than 7.
\end{thm}

\begin{thm}\label{Rook'sgraph}
The twin-width of a $(n \times m)$ Rook's graph is less than equal to $2(m - 1)$.

\end{thm}

\subsection{Organization}
In Section 2, we introduce our necessary definitions, notations and abreviations. In Section 3, we survey the known results of twin-width of finite graphs. We study the behaviour of twin-width under graph operations in Section 4. In Section 5, we compute the twin-width of finite graphs with 4 and 5 vertices. In Section 6, we provide upper bounds of twin-width of King's graph and Rook's graph.

\section{Preliminaries: some definitions, notations and abreviations:}

A \textit{trigraph} $G$ is a graph with a vertex set $V(G)$, a black edge set $E(G)$, and a red edge set $R(G)$ (the error
edges), where $E(G)$ and $R(G)$ are disjoint. The set of neighbours of a vertex $v$ in a
trigraph $G$, denoted by $N_G(v)$, consists of all the vertices adjacent to $v$ by a black or red edge. The degree of a vertex 
$v$ is defined by the number $\mid N_G(v)\mid$. A $d$-trigraph is a
trigraph $G$ such that the red graph $(V(G), R(G))$ has degree at most $d$. In this situation, we also
say that the trigraph has red degree at most $d$.

A \textit{contraction} or \textit{identification} in a
trigraph $G$ consists of merging two (non-necessarily adjacent) vertices $u$ and $v$ into a single
vertex $w$, and defining the edges of $G'$ (the new graph after contraction) in the following way: Every vertex of the symmetric
difference $N_G(u)\bigtriangleup N_G(v)$ is linked to $w$ by a red edge. Every vertex $x$ of the intersection $N_G(u)\cap N_G(v)$ is linked to $w$ by a black edge if both $ux\in E(G)$ and $vx\in E(G)$, and by a red edge otherwise. The rest of the edges (not incident to $u$ or $v$) remain unchanged.  Also, the vertices $u$ and $v$ (together with the edges incident to these vertices) are removed from the trigraph. 

A \textit{sequence of $d$-contractions} is a sequence of $d$-trigraphs $G_n, G_{n-1},\cdots, G_1$, where $G_n = G$, $G_1 = K_1$, where $K_1$ is the graph on a single vertex, and $G_{i-1}$ is obtained from $G_i$ by performing a single contraction of two (non-necessarily adjacent) vertices. We observe that $G_i$ has precisely $i$ vertices, for every $i\in \{1,\cdots, n\}$. The twin-width of $G$, denoted by $tww(G)$, is the minimum integer $d$ such that $G$ admits a $d$-sequence.

Now, we provide an example of a sequence contractions of a finite graph. In the sequence of graphs depicted below, we start with the given finite graph in the extreme left end and we label the vertices by $a, b, c, d, e, f, g$. The next diagram is the result of the contraction of the vertices $e$ and $f$ and in the resulting graph we label the new vertex by $ef$. In this way, we obtain a sequence of graphs by gradually contacting other vertices. The graph in the extreme left end of the second line of this sequence is obtained by contracting the vertices $ad$ and $g$ in the graph depicted in the extreme right end of the first line of the sequence.

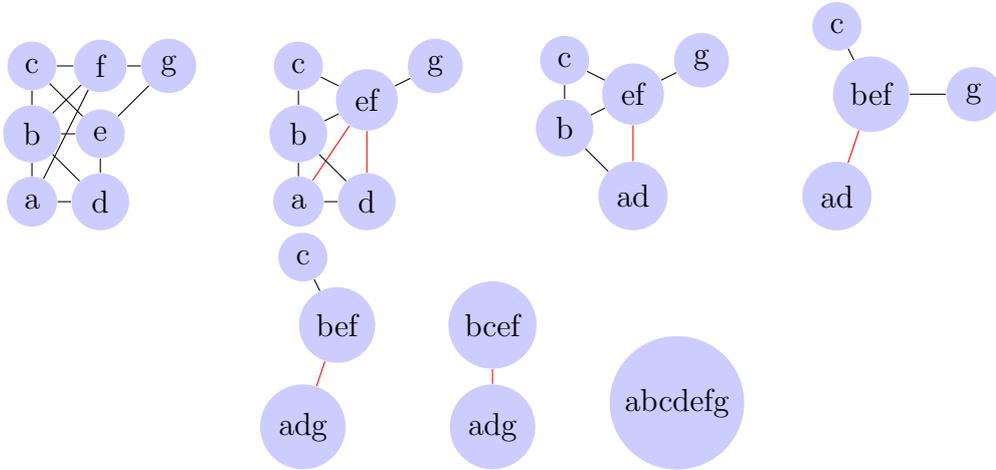
\begin{center}
\begin{tikzpicture}  
  [scale=.9,auto=center,every node/.style={circle,fill=blue!20}] 
    
  \node (a1) at (0,0) {a};  
  \node (a2) at (0,1)  {b};  
  \node (a3) at (0,2)  {c};  
  \node (a4) at (1,0) {d};  
  \node (a5) at (1,1)  {e};  
  \node (a6) at (1,2)  {f};  
  \node (a7) at (2,2)  {g};  
  
  \draw (a1) -- (a2); 
  \draw (a1) -- (a4);
  \draw (a1) -- (a6);
  \draw (a2) -- (a3); 
  \draw (a2) -- (a4); 
  \draw (a2) -- (a5);
  \draw (a2) -- (a6);
 \draw (a3) -- (a5);
 \draw (a3) -- (a6);
 \draw (a4) -- (a5);  
 \draw (a5) -- (a7);
 \draw (a6) -- (a7);  
\end{tikzpicture} 
\qquad
\begin{tikzpicture}  
  [scale=.9,auto=center,every node/.style={circle,fill=blue!20}] 
    
  \node (a1) at (0,0) {a};  
  \node (a2) at (0,1)  {b};  
  \node (a3) at (0,2)  {c};  
  \node (a4) at (1,0) {d};  
  \node (a5) at (1,1.5)  {ef};  
  \node (a7) at (2,2)  {g};  
  
  \draw (a1) -- (a2); 
  \draw (a1) -- (a4);
  \draw[red] (a1) -- (a5);
  \draw (a2) -- (a3); 
  \draw (a2) -- (a4); 
  \draw (a2) -- (a5);
  \draw (a3) -- (a5);

 \draw[red] (a4) -- (a5);  
 \draw (a5) -- (a7);

\end{tikzpicture}
\qquad
\begin{tikzpicture}  
  [scale=.9,auto=center,every node/.style={circle,fill=blue!20}] 
    
  \node (a1) at (1,0) {ad};  
  \node (a2) at (0,1)  {b};  
  \node (a3) at (0,2)  {c};  
 
  \node (a5) at (1,1.5)  {ef};  
  \node (a7) at (2,2)  {g};  
  
  \draw (a1) -- (a2); 
  \draw[red] (a1) -- (a5);
  \draw (a2) -- (a3); 

  \draw (a2) -- (a5);
  \draw (a3) -- (a5);

 \draw (a5) -- (a7);

\end{tikzpicture}
\qquad
\begin{tikzpicture}  
  [scale=.9,auto=center,every node/.style={circle,fill=blue!20}] 
    
  \node (a1) at (0,-0.5) {ad};  
  \node (a2) at (0.5,1)  {bef};  
  \node (a3) at (0,2)  {c};  
  \node (a7) at (2,1)  {g};  
  
  \draw[red] (a1) -- (a2); 
  \draw (a2) -- (a3); 
  \draw (a2) -- (a7);

\end{tikzpicture}
\qquad
\begin{tikzpicture}  
  [scale=.9,auto=center,every node/.style={circle,fill=blue!20}] 
    
  \node (a1) at (0,-0.5) {adg};  
  \node (a2) at (0.5,1)  {bef};  
  \node (a3) at (0,2)  {c};

  \draw[red] (a1) -- (a2); 
  \draw (a2) -- (a3);

\end{tikzpicture}
\qquad
\begin{tikzpicture}  
  [scale=.9,auto=center,every node/.style={circle,fill=blue!20}] 
    
  \node (a1) at (0,-0.5) {adg};  
  \node (a2) at (0,1)  {bcef};

  \draw[red] (a1) -- (a2); 
\end{tikzpicture}
\qquad
\begin{tikzpicture}  
  [scale=.4,auto=center,every node/.style={circle,fill=blue!20}] 
    
  \node (a1) at (0,0) {abcdefg};

\end{tikzpicture}

\end{center}
Moreover, we will draw every contraction sequence by this fashion in this article. We end this section 
by defining twin-width of an infinite graph. It is defined by the maximum of the twin-widths of its induced finite subgraphs. 

\section{A survey of the known results of twin-width of finite graphs}

In this section, we survey the results regarding the twin-width of complete graphs, planar graphs, graphs with at most one cycle, Caterpillar tree and Paley graph. 

\begin{thm}
The complete graph with $n$-vertices, denoted by $K_n$,  and the complete bipartite with ${n,m}$-vertices, denoted by $K_{n,m}$,  have twin-width zero.
\end{thm}

\begin{thm}\cite{H22}
The twin-width of any simple planar graph $G$ is at most 9.
\end{thm}

\begin{thm}\label{onecycle}\cite{AHKO22}
If every component of a graph $G$ has at most one cycle,
then $tww(G)\leq 2$.
\end{thm}

We obtain the following corollary from the above-mentioned theorem.
\begin{cor}\label{cyclicgraph}
The cyclic graph with $n$-vertices (denoted by $C_n$) has twin-width less than equal to 2.
\end{cor}

A \textit{caterpillar tree} is a tree in which all the vertices are within distance 1 of a central path.
We draw an example of a Caterpillar tree below. 

\begin{center}
\begin{tikzpicture} 
[scale=1.5,auto=center,every node/.style={circle,fill=blue!20}] 
\node (a1) at (0,0) {a}; 
\node (a2) at (1,0) {b}; 
\node (a3) at (2,0){c};
\node (a4) at (2,-1) {d};
\node (a5) at (3,0) {e};
\node (a6) at (2.5,1) {f};
\node (a7) at (3,1) {g};
\node (a8) at (3.5,1) {h};
\node (a9) at (2.5,-1) {i};
\node (a10) at (3.5,-1){j};
\node (a11) at (4,0) {k};

\draw (a1) -- (a2);
\draw (a2) -- (a3);
\draw (a3) -- (a4);
\draw (a3) -- (a5);
\draw (a5) -- (a6);
\draw (a5) -- (a7);
\draw (a5) -- (a8);
\draw (a5) -- (a9);
\draw (a5) -- (a10);
\draw (a5) -- (a11);
\end{tikzpicture}
\end{center}

\begin{thm}\label{caterpillar}\cite{AHKO22}
For a tree $T$, $tww(T)\leq 1$ if and only if $T$ is a Caterpillar tree. 
\end{thm}
The above-mentioned theorem gives rise to the following corollary. 
\begin{cor}\label{pathgraph}
The path graph with $n$ vertices, denoted by $P_n$, has twin-width less than equal to 1. 
\end{cor}

Let $q$ be a prime power such that $q\equiv 1$ (mod 4). The Paley graph of order $q$, denoted by $P(q)$, is defined as follows: The
vertices of the graph are the elements of the field $\textbf{F}_q$ and the vertices $i$ and $j$ are adjacent if $(j-i)$ is a quadratic residue in $\textbf{F}_q$. 
%We draw below a Paley graph of order 13.
%\begin{figure}[!ht]
 % \centering
 % \includegraphics[width=0.5\textwidth,natwidth=610,natheight=642]{Paley13.jpg}
%\end{figure}

\begin{thm}\cite{AHKO22}
For each prime $q$ with  $q\equiv 1$ (mod 4), the Paley graph
$P(q)$ has twin-width exactly $(q-1)/2$. 
\end{thm}

\begin{rem}
The twin-width of the class of Paley graphs is unbounded.
\end{rem}

Finally, we end this section by mentioning results on groups.  We say that a finitely generated group has bounded or unbounded twin-width if one of its Cayley graphs have bounded or unbounded twin-width. 

\begin{thm}\cite{BGTT22}
The solvable, hyperbolic, ordered finitely generated groups have finite twin-width.
\end{thm}

\begin{thm}\cite{BGTT22}
There is a finitely generated group with infinite twin-width.
\end{thm}

\section{Graph operations and their twin-width}

In this section, we study the behaviour of twin-width under graph operations, like complement graph, dual graph and line graph.

 \subsection{Complement of a graph}\label{subsection4.3}
 The \textit{complement of a graph} $G$ is a graph $H$ on the same vertices such that two distinct vertices of $H$ are adjacent if and only if they are not adjacent in $G$.  We obtain the following theorem from \cite{BKTW20} (see Subsection 4.1).
 
 \begin{thm}\label{complementgraph}
 Twin-width is invariant under complementation. 
 \end{thm}

\subsection{Dual Graph}\label{subsection4.1}

The \textit{dual graph} of a planar graph $G$ is a graph that has a vertex for each face of $G$. The dual graph has an edge for each pair of faces in $G$ that are separated from each other by an edge, and a self-loop when the same face appears on both sides of an edge. Now, we  
prove Theorem \ref{dual}.

\vspace{5mm}

\textbf{Proof of Theorem \ref{dual}:}
Let $G$ be the following graph: 
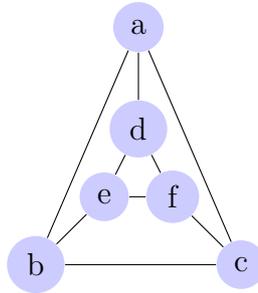
\begin{center}
\begin{tikzpicture} 
[scale=.9,auto=center,every node/.style={circle,fill=blue!20}] 
\node (a1) at (-0.5,-0.5) {e}; 
\node (a2) at (0.5,-0.5) {f}; 
\node (a3) at (0,0.5){d};
\node (a4) at (-1.5,-1.5) {b};
\node (a5) at (1.5,-1.5) {c};
\node (a6) at (0,2) {a};

\draw (a1) -- (a2);
\draw (a2) -- (a3);
\draw (a3) -- (a1);
\draw (a4) -- (a5);
\draw (a5) -- (a6);
\draw (a6) -- (a4);
\draw (a1) -- (a4);
\draw (a2) -- (a5);
\draw (a3) -- (a6);
\end{tikzpicture}
\end{center}
If we contract any two vertices of $G$, it generates a red edge. 
Therefore, the twin-width of $G$ is greater than equal to 1. Now, we
compute the dual graph of $G$. We denote the triangular region `def' by $r1$, the trapezium `bcfe' by $r2$, the trapezium `adfc' by $r3$, the trapezium
`abed' by $r4$ and the region outside $\triangle abc$ by $r5$. Therefore, the dual of $G$, denoted by $G^*$, will be the following graph: 
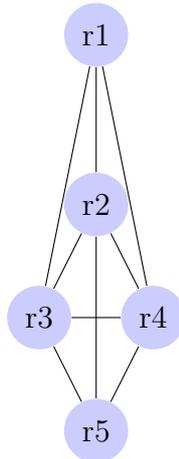
\begin{center} 
\begin{tikzpicture} 
[scale=1.5,auto=center,every node/.style={circle,fill=blue!20}] 
\node (a1) at (-0.5,-0.5) {r3}; 
\node (a2) at (0.5,-0.5) {r4}; 
\node (a3) at (0,0.5){r2};
\node (a4) at (0,-1.5) {r5};
\node (a5) at (0,2) {r1};

\draw (a1) -- (a2);
\draw (a2) -- (a3);
\draw (a3) -- (a1);
\draw (a4) -- (a1);
\draw (a4) -- (a2);
\draw (a4) -- (a3);
\draw (a5) -- (a1);
\draw (a5) -- (a2);
\draw (a5) -- (a3);
\end{tikzpicture}
\end{center}

We apply the following $0$-contraction sequence to $G^*$. 

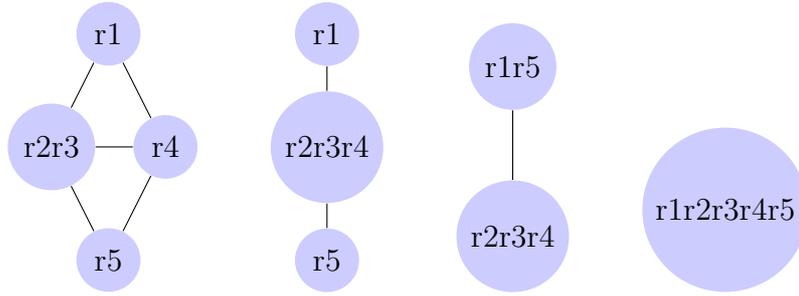
\begin{center} 
\begin{tikzpicture} 
[scale=1.5,auto=center,every node/.style={circle,fill=blue!20}] 
\node (a1) at (-0.5,-0.5) {r2r3}; 
\node (a2) at (0.5,-0.5) {r4}; 
%\node (a3) at (0,0.5){r2};
\node (a4) at (0,-1.5) {r5};
\node (a5) at (0,0.5) {r1};

\draw (a1) -- (a2);
%\draw (a2) -- (a3);
%\draw (a3) -- (a1);
\draw (a4) -- (a1);
\draw (a4) -- (a2);
%\draw (a4) -- (a3);
\draw (a5) -- (a1);
\draw (a5) -- (a2);
%\draw (a5) -- (a3);
\end{tikzpicture}
\qquad
\begin{tikzpicture} 
[scale=1.5,auto=center,every node/.style={circle,fill=blue!20}] 
\node (a1) at (0,0) {r2r3r4}; 
%\node (a2) at (0.5,-0.5) {r4}; 
%\node (a3) at (0,0.5){r2};
\node (a4) at (0,-1) {r5};
\node (a5) at (0,1) {r1};

\draw (a1) -- (a4);
%\draw (a2) -- (a3);
%\draw (a3) -- (a1);
\draw (a5) -- (a1);
%\draw (a4) -- (a2);
%\draw (a4) -- (a3);
%\draw (a5) -- (a1);
%\draw (a5) -- (a2);
%\draw (a5) -- (a3);
\end{tikzpicture}
\qquad
\begin{tikzpicture} 
[scale=1.5,auto=center,every node/.style={circle,fill=blue!20}] 
\node (a1) at (0,0) {r2r3r4}; 

\node (a4) at (0,1.5) {r1r5};

\draw (a1) -- (a4);

\end{tikzpicture}
\qquad
\begin{tikzpicture} 
[scale=1.5,auto=center,every node/.style={circle,fill=blue!20}] 
\node (a1) at (0,0) {r1r2r3r4r5}; 
%\node (a4) at (0,2) {r1r5};
\end{tikzpicture}
\end{center}
Therefore, $G^*$ has twin-width zero. Hence, the twin-widths of $G$ and $G^*$ are 
different.
\hfill\(\Box\)

\subsection{Line graph of a graph}\label{subsection4.2}

 The \textit{line graph} of an undirected graph $G$ is another graph $L(G)$ that represents the adjacencies between edges of $G$. $L(G)$ is constructed in the following way: for each edge in $G$, make a vertex in $L(G)$; for every two edges in $G$ that have a vertex in common, make an edge between their  corresponding vertices in $L(G)$. Now, we prove the Theorem \ref{linegraph}.
 
 \vspace{5mm}

\textbf{Proof of Theorem \ref{linegraph}:} Let, $G$ be the following graph:
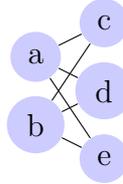
\begin{center}
\begin{tikzpicture}  
  [scale=.9,auto=center,every node/.style={circle,fill=blue!20}] 
    
  \node (a1) at (0,1.5) {a};  
  \node (a2) at (0,0.5)  {b};  
  \node (a3) at (1,2)  {c};  
  \node (a4) at (1,1) {d};  
  \node (a5) at (1,0)  {e};

  \draw (a1) -- (a3); 
  \draw (a1) -- (a4);
  \draw (a1) -- (a5);
  \draw (a2) -- (a3); 
  \draw (a2) -- (a4); 
  \draw (a2) -- (a5);
\end{tikzpicture} 
\end{center} 
Since $G$ is a complete bipartite graph, it has twin-width zero. However, the line graph of $G$, denoted by $L(G)$, is the following graph:
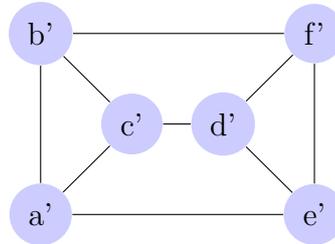
\begin{center}
\begin{tikzpicture}  
  [scale=1.2,auto=center,every node/.style={circle,fill=blue!20}] 
    
  \node (a1) at (0,0) {a'};  
  \node (a2) at (0,2)  {b'};  
  \node (a3) at (1,1)  {c'};  
  \node (a4) at (2,1) {d'};  
  \node (a5) at (3,0)  {e'};  
  \node (a6) at (3,2) {f'};

  \draw (a1) -- (a2); 
  \draw (a1) -- (a3);
  \draw (a1) -- (a5);
  \draw (a2) -- (a3); 
  \draw (a2) -- (a6); 
  \draw (a3) -- (a4);
  \draw (a4) -- (a5);
  \draw (a4) -- (a6);
  \draw (a5) -- (a6);
\end{tikzpicture}  
\end{center}
We observe that if we contract any two vertices of $L(G)$, it generates a red edge. Therefore, the twin-width of $L(G)$ is greater than equal to $1$, which implies that the twin-width is not preserved under taking line graph. \hfill\(\Box\)

\section{Computation of twin-width of finite graphs with 4 and 5 vertices}

In this section, we prove Theorem \ref{4vertices} and Theorem \ref{5vertices}.

\vspace{5mm}

\textbf{Proof of Theorem \ref{4vertices}:} First, we make a list of graphs which are disconnected.
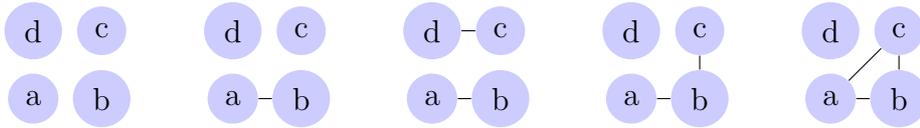
\begin{center}
\begin{tikzpicture} 
[scale=.9,auto=center,every node/.style={circle,fill=blue!20}] 
\node (a1) at (0,0) {a}; 
\node (a2) at (1,0) {b}; 
\node (a3) at (1,1) {c}; 
\node (a4) at (0,1) {d};
\end{tikzpicture}
\qquad
\begin{tikzpicture} 
[scale=.9,auto=center,every node/.style={circle,fill=blue!20}] 
\node (a1) at (0,0) {a}; 
\node (a2) at (1,0) {b}; 
\node (a3) at (1,1) {c}; 
\node (a4) at (0,1) {d};

\draw (a1) -- (a2);

\end{tikzpicture}
\qquad
\begin{tikzpicture} 
[scale=.9,auto=center,every node/.style={circle,fill=blue!20}] 
\node (a1) at (0,0) {a}; 
\node (a2) at (1,0) {b}; 
\node (a3) at (1,1) {c}; 
\node (a4) at (0,1) {d};

\draw (a1) -- (a2);
\draw (a4) -- (a3);

\end{tikzpicture}
\qquad
\begin{tikzpicture} 
[scale=.9,auto=center,every node/.style={circle,fill=blue!20}] 
\node (a1) at (0,0) {a}; 
\node (a2) at (1,0) {b}; 
\node (a3) at (1,1) {c}; 
\node (a4) at (0,1) {d};

\draw (a1) -- (a2);

\draw (a2) -- (a3);

\end{tikzpicture}
\qquad
\begin{tikzpicture} 
[scale=.9,auto=center,every node/.style={circle,fill=blue!20}] 
\node (a1) at (0,0) {a}; 
\node (a2) at (1,0) {b}; 
\node (a3) at (1,1) {c}; 
\node (a4) at (0,1) {d};

\draw (a1) -- (a2);
\draw (a1) -- (a3);

\draw (a2) -- (a3);

\end{tikzpicture}
\end{center}
Since the twin-width of a (disconnected) graph is the maximum of the twin-widths of its components, it is easy to see that the twin-widths of the disconnected graphs with 4 vertices are zero. Now, we make a list of graphs which are Caterpillar tree. 
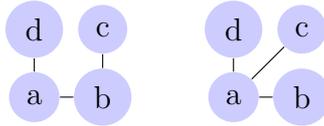
\begin{center}
\begin{tikzpicture} 
[scale=.9,auto=center,every node/.style={circle,fill=blue!20}] 
\node (a1) at (0,0) {a}; 
\node (a2) at (1,0) {b}; 
\node (a3) at (1,1) {c}; 
\node (a4) at (0,1) {d};

\draw (a1) -- (a2);
\draw (a1) -- (a4);
\draw (a2) -- (a3);
\end{tikzpicture}
\qquad
\begin{tikzpicture} 
[scale=.9,auto=center,every node/.style={circle,fill=blue!20}] 
\node (a1) at (0,0) {a}; 
\node (a2) at (1,0) {b}; 
\node (a3) at (1,1) {c}; 
\node (a4) at (0,1) {d};

\draw (a1) -- (a2);
\draw (a1) -- (a4);
\draw (a1) -- (a3);
\end{tikzpicture}
\end{center}
By Theorem \ref{caterpillar}, these two graphs have twin-width less than equal to 1. Finally, we make a list of graphs whose complement graphs are disconnected. 
 
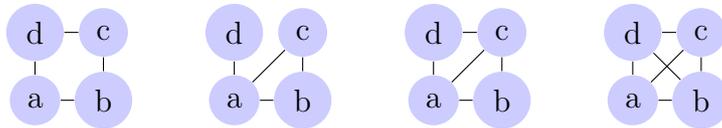
\begin{center}
\begin{tikzpicture} 
[scale=.9,auto=center,every node/.style={circle,fill=blue!20}] 
\node (a1) at (0,0) {a}; 
\node (a2) at (1,0) {b}; 
\node (a3) at (1,1) {c}; 
\node (a4) at (0,1) {d};

\draw (a1) -- (a2);
\draw (a1) -- (a4);
\draw (a2) -- (a3);
\draw (a3) -- (a4);
\end{tikzpicture}
\qquad
\begin{tikzpicture} 
[scale=.9,auto=center,every node/.style={circle,fill=blue!20}] 
\node (a1) at (0,0) {a}; 
\node (a2) at (1,0) {b}; 
\node (a3) at (1,1) {c}; 
\node (a4) at (0,1) {d};

\draw (a1) -- (a2);
\draw (a1) -- (a3);
\draw (a1) -- (a4);
\draw (a2) -- (a3);
\end{tikzpicture}
\qquad
\begin{tikzpicture} 
[scale=.9,auto=center,every node/.style={circle,fill=blue!20}] 
\node (a1) at (0,0) {a}; 
\node (a2) at (1,0) {b}; 
\node (a3) at (1,1) {c}; 
\node (a4) at (0,1) {d};

\draw (a1) -- (a2);
\draw (a1) -- (a3);
\draw (a1) -- (a4);
\draw (a2) -- (a3);
\draw (a3) -- (a4);
\end{tikzpicture}
\qquad
\begin{tikzpicture} 
[scale=.9,auto=center,every node/.style={circle,fill=blue!20}] 
\node (a1) at (0,0) {a}; 
\node (a2) at (1,0) {b}; 
\node (a3) at (1,1) {c}; 
\node (a4) at (0,1) {d};

\draw (a1) -- (a2);
\draw (a1) -- (a3);
\draw (a1) -- (a4);
\draw (a2) -- (a3);
\draw (a2) -- (a4);
\draw (a3) -- (a4);
\end{tikzpicture}
\end{center}
Since the twin-width of a graph is same as the twin-width of the complement graph by Theorem \ref{complementgraph} and the twin-width of a disconnected graph with 4 vertices is less than equal to 1,
we obtain that the twin-widths of the graphs in the above list is zero. Hence, we have our theorem.
\hfill\(\Box\)

\vspace{5mm}

\textbf{Proof of Theorem \ref{5vertices}:} First, we make a list of  disconnected simple graphs with 5 vertices:

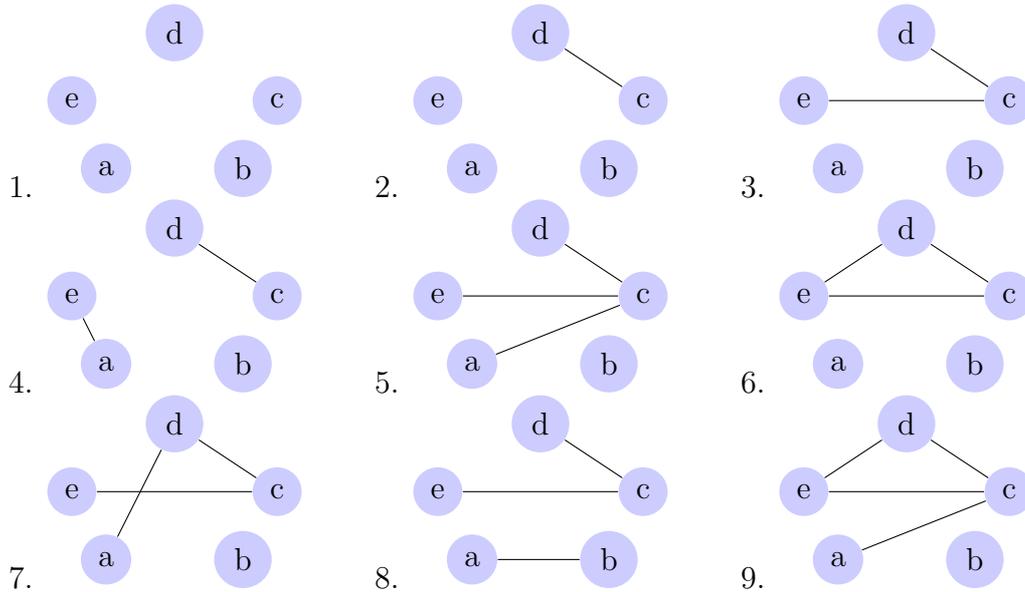
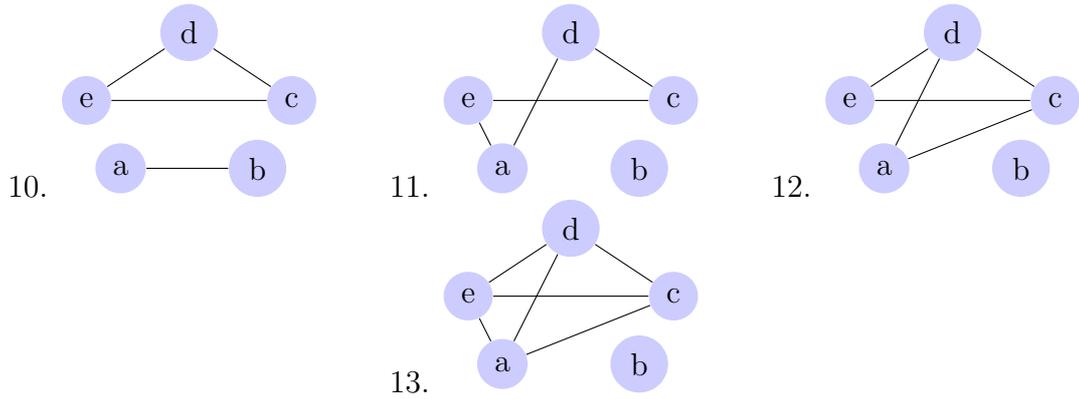
\begin{center}
1.
\begin{tikzpicture} 
[scale=.9,auto=center,every node/.style={circle,fill=blue!20}] 
\node (a1) at (-1,0) {a}; 
\node (a2) at (1,0) {b}; 
\node (a3) at (3/2,1) {c}; 
\node (a4) at (0,2) {d};
\node (a5) at (-3/2,1) {e};  
\end{tikzpicture}
\qquad
2.
\begin{tikzpicture} 
[scale=.9,auto=center,every node/.style={circle,fill=blue!20}] 
\node (a1) at (-1,0) {a}; 
\node (a2) at (1,0) {b}; 
\node (a3) at (3/2,1) {c}; 
\node (a4) at (0,2) {d};
\node (a5) at (-3/2,1) {e};

\draw (a3) -- (a4);

\end{tikzpicture}
\qquad
3.
\begin{tikzpicture} 
[scale=.9,auto=center,every node/.style={circle,fill=blue!20}] 
\node (a1) at (-1,0) {a}; 
\node (a2) at (1,0) {b}; 
\node (a3) at (3/2,1) {c}; 
\node (a4) at (0,2) {d};
\node (a5) at (-3/2,1) {e};

\draw (a3) -- (a4);
\draw (a3) -- (a5);

\end{tikzpicture}
\qquad
4.
\begin{tikzpicture} 
[scale=.9,auto=center,every node/.style={circle,fill=blue!20}] 
\node (a1) at (-1,0) {a}; 
\node (a2) at (1,0) {b}; 
\node (a3) at (3/2,1) {c}; 
\node (a4) at (0,2) {d};
\node (a5) at (-3/2,1) {e};

\draw (a1) -- (a5);

\draw (a3) -- (a4);
\end{tikzpicture}
\qquad
5.
\begin{tikzpicture} 
[scale=.9,auto=center,every node/.style={circle,fill=blue!20}] 
\node (a1) at (-1,0) {a}; 
\node (a2) at (1,0) {b}; 
\node (a3) at (3/2,1) {c}; 
\node (a4) at (0,2) {d};
\node (a5) at (-3/2,1) {e};

\draw (a1) -- (a3);

\draw (a3) -- (a4);
\draw (a3) -- (a5);

\end{tikzpicture}
\qquad
6.
\begin{tikzpicture} 
[scale=.9,auto=center,every node/.style={circle,fill=blue!20}] 
\node (a1) at (-1,0) {a}; 
\node (a2) at (1,0) {b}; 
\node (a3) at (3/2,1) {c}; 
\node (a4) at (0,2) {d};
\node (a5) at (-3/2,1) {e};

\draw (a3) -- (a4);
\draw (a3) -- (a5);
\draw (a4) -- (a5);

\end{tikzpicture}
\qquad
7.
\begin{tikzpicture} 
[scale=.9,auto=center,every node/.style={circle,fill=blue!20}] 
\node (a1) at (-1,0) {a}; 
\node (a2) at (1,0) {b}; 
\node (a3) at (3/2,1) {c}; 
\node (a4) at (0,2) {d};
\node (a5) at (-3/2,1) {e};

\draw (a1) -- (a4);

\draw (a3) -- (a4);
\draw (a3) -- (a5);

\end{tikzpicture}
\qquad
8.
\begin{tikzpicture} 
[scale=.9,auto=center,every node/.style={circle,fill=blue!20}] 
\node (a1) at (-1,0) {a}; 
\node (a2) at (1,0) {b}; 
\node (a3) at (3/2,1) {c}; 
\node (a4) at (0,2) {d};
\node (a5) at (-3/2,1) {e};  

\draw (a1) -- (a2);

\draw (a3) -- (a4);
\draw (a3) -- (a5);
%\draw (a4) -- (a5);

\end{tikzpicture}
\qquad
9.
\begin{tikzpicture} 
[scale=.9,auto=center,every node/.style={circle,fill=blue!20}] 
\node (a1) at (-1,0) {a}; 
\node (a2) at (1,0) {b}; 
\node (a3) at (3/2,1) {c}; 
\node (a4) at (0,2) {d};
\node (a5) at (-3/2,1) {e};  

%\draw (a1) -- (a2);
\draw (a1) -- (a3);
\draw (a3) -- (a4);
\draw (a3) -- (a5);
\draw (a4) -- (a5);
\end{tikzpicture}
\qquad
10.
\begin{tikzpicture} 
[scale=.9,auto=center,every node/.style={circle,fill=blue!20}] 
\node (a1) at (-1,0) {a}; 
\node (a2) at (1,0) {b}; 
\node (a3) at (3/2,1) {c}; 
\node (a4) at (0,2) {d};
\node (a5) at (-3/2,1) {e};  

\draw (a1) -- (a2);

\draw (a3) -- (a4);
\draw (a3) -- (a5);
\draw (a4) -- (a5);

\end{tikzpicture}
\qquad
11.
\begin{tikzpicture} 
[scale=.9,auto=center,every node/.style={circle,fill=blue!20}] 
\node (a1) at (-1,0) {a}; 
\node (a2) at (1,0) {b}; 
\node (a3) at (3/2,1) {c}; 
\node (a4) at (0,2) {d};
\node (a5) at (-3/2,1) {e};  

%\draw (a1) -- (a2);
%\draw (a1) -- (a3);
\draw (a1) -- (a4);
\draw (a1) -- (a5);
%\draw (a2) -- (a3);

\draw (a3) -- (a4);
\draw (a3) -- (a5);
%\draw (a4) -- (a5);

\end{tikzpicture}
\qquad
12.
\begin{tikzpicture} 
[scale=.9,auto=center,every node/.style={circle,fill=blue!20}] 
\node (a1) at (-1,0) {a}; 
\node (a2) at (1,0) {b}; 
\node (a3) at (3/2,1) {c}; 
\node (a4) at (0,2) {d};
\node (a5) at (-3/2,1) {e};  

%\draw (a1) -- (a2);
\draw (a1) -- (a3);
\draw (a1) -- (a4);
%\draw (a1) -- (a5);

\draw (a3) -- (a4);
\draw (a3) -- (a5);
\draw (a4) -- (a5);

\end{tikzpicture}
\qquad
13.
\begin{tikzpicture} 
[scale=.9,auto=center,every node/.style={circle,fill=blue!20}] 
\node (a1) at (-1,0) {a}; 
\node (a2) at (1,0) {b}; 
\node (a3) at (3/2,1) {c}; 
\node (a4) at (0,2) {d};
\node (a5) at (-3/2,1) {e};  

%\draw (a1) -- (a2);
\draw (a1) -- (a3);
\draw (a1) -- (a4);
\draw (a1) -- (a5);

\draw (a3) -- (a4);
\draw (a3) -- (a5);
\draw (a4) -- (a5);

\end{tikzpicture}
\end{center}
Since the twin-width of a disconnected graph is the maximum of the twin-widths of its components and the twin-width of a simple graph with 4 vertices is less than equal to 1 by Theorem \ref{4vertices},   the twin-width of the disconnected graphs with 5 vertices are less than equal to 1. Now, we make a  list of graphs which are Caterpillar tree and have 5 vertices. 

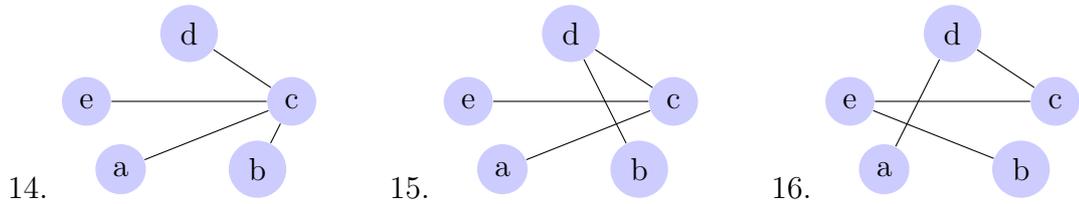
\begin{center}
14.
\begin{tikzpicture} 
[scale=.9,auto=center,every node/.style={circle,fill=blue!20}] 
\node (a1) at (-1,0) {a}; 
\node (a2) at (1,0) {b}; 
\node (a3) at (3/2,1) {c}; 
\node (a4) at (0,2) {d};
\node (a5) at (-3/2,1) {e};  

%\draw (a1) -- (a2);
\draw (a1) -- (a3);

\draw (a2) -- (a3);

\draw (a3) -- (a4);
\draw (a3) -- (a5);
\end{tikzpicture}
\qquad
15.
\begin{tikzpicture} 
[scale=.9,auto=center,every node/.style={circle,fill=blue!20}] 
\node (a1) at (-1,0) {a}; 
\node (a2) at (1,0) {b}; 
\node (a3) at (3/2,1) {c}; 
\node (a4) at (0,2) {d};
\node (a5) at (-3/2,1) {e};  

%\draw (a1) -- (a2);
\draw (a1) -- (a3);

%\draw (a2) -- (a3);
\draw (a2) -- (a4);
%\draw (a2) -- (a5);
\draw (a3) -- (a4);
\draw (a3) -- (a5);
%\draw (a4) -- (a5);

\end{tikzpicture}
\qquad
16.
\begin{tikzpicture} 
[scale=.9,auto=center,every node/.style={circle,fill=blue!20}] 
\node (a1) at (-1,0) {a}; 
\node (a2) at (1,0) {b}; 
\node (a3) at (3/2,1) {c}; 
\node (a4) at (0,2) {d};
\node (a5) at (-3/2,1) {e};  

%\draw (a1) -- (a2);
%\draw (a1) -- (a3);
\draw (a1) -- (a4);
%\draw (a1) -- (a5);
%\draw (a2) -- (a3);
%\draw (a2) -- (a4);
\draw (a2) -- (a5);
\draw (a3) -- (a4);
\draw (a3) -- (a5);
%\draw (a4) -- (a5);

\end{tikzpicture}
\end{center}
Since a Caterpillar tree has twin-width less than equal to 1 by \ref{caterpillar}, the graphs in the above list have twin-width less than equal to 1. Next, we make a list of graphs which have at most one cycle. 

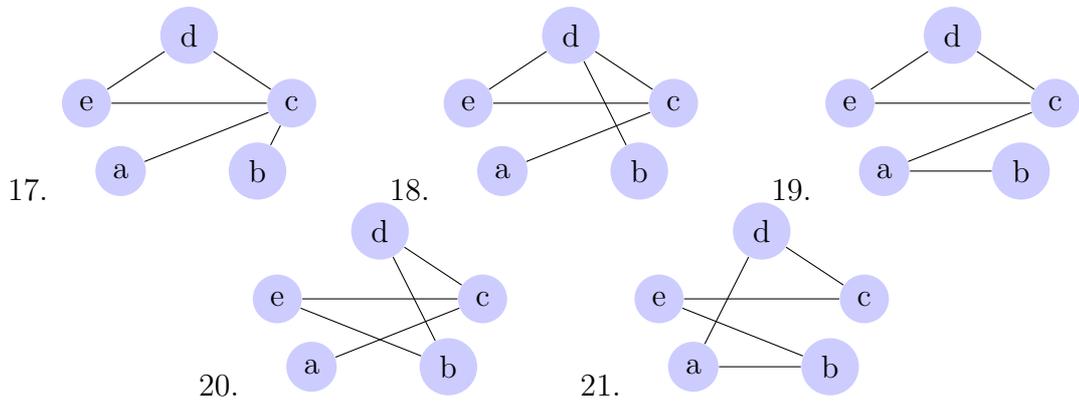
\begin{center}
17.
\begin{tikzpicture} 
[scale=.9,auto=center,every node/.style={circle,fill=blue!20}] 
\node (a1) at (-1,0) {a}; 
\node (a2) at (1,0) {b}; 
\node (a3) at (3/2,1) {c}; 
\node (a4) at (0,2) {d};
\node (a5) at (-3/2,1) {e};  

%\draw (a1) -- (a2);
\draw (a1) -- (a3);

\draw (a2) -- (a3);

\draw (a3) -- (a4);
\draw (a3) -- (a5);
\draw (a4) -- (a5);

\end{tikzpicture}
\qquad
18.
\begin{tikzpicture} 
[scale=.9,auto=center,every node/.style={circle,fill=blue!20}] 
\node (a1) at (-1,0) {a}; 
\node (a2) at (1,0) {b}; 
\node (a3) at (3/2,1) {c}; 
\node (a4) at (0,2) {d};
\node (a5) at (-3/2,1) {e};  

%\draw (a1) -- (a2);
\draw (a1) -- (a3);
%\draw (a1) -- (a4);
%\draw (a1) -- (a5);
%\draw (a2) -- (a3);
\draw (a2) -- (a4);
%\draw (a2) -- (a5);
\draw (a3) -- (a4);
\draw (a3) -- (a5);
\draw (a4) -- (a5);

\end{tikzpicture}
\qquad
19.
\begin{tikzpicture} 
[scale=.9,auto=center,every node/.style={circle,fill=blue!20}] 
\node (a1) at (-1,0) {a}; 
\node (a2) at (1,0) {b}; 
\node (a3) at (3/2,1) {c}; 
\node (a4) at (0,2) {d};
\node (a5) at (-3/2,1) {e};  

\draw (a1) -- (a2);
\draw (a1) -- (a3);

\draw (a3) -- (a4);
\draw (a3) -- (a5);
\draw (a4) -- (a5);

\end{tikzpicture}
\qquad
20.
\begin{tikzpicture} 
[scale=.9,auto=center,every node/.style={circle,fill=blue!20}] 
\node (a1) at (-1,0) {a}; 
\node (a2) at (1,0) {b}; 
\node (a3) at (3/2,1) {c}; 
\node (a4) at (0,2) {d};
\node (a5) at (-3/2,1) {e};  

%\draw (a1) -- (a2);
\draw (a1) -- (a3);
%\draw (a1) -- (a4);
%\draw (a1) -- (a5);
%\draw (a2) -- (a3);
\draw (a2) -- (a4);
\draw (a2) -- (a5);
\draw (a3) -- (a4);
\draw (a3) -- (a5);
%\draw (a4) -- (a5);

\end{tikzpicture}
\qquad
21.
\begin{tikzpicture} 
[scale=.9,auto=center,every node/.style={circle,fill=blue!20}] 
\node (a1) at (-1,0) {a}; 
\node (a2) at (1,0) {b}; 
\node (a3) at (3/2,1) {c}; 
\node (a4) at (0,2) {d};
\node (a5) at (-3/2,1) {e};  

\draw (a1) -- (a2);
%\draw (a1) -- (a3);
\draw (a1) -- (a4);
%\draw (a1) -- (a5);
%\draw (a2) -- (a3);
%\draw (a2) -- (a4);
\draw (a2) -- (a5);
\draw (a3) -- (a4);
\draw (a3) -- (a5);
%\draw (a4) -- (a5);

\end{tikzpicture}
\end{center}
Since the graphs with at most one cycle have twin-width less than equal to 2, the graphs in the above-mentioned list have twin-width less than equal to 2. Now, we make a list of graphs which can be reduced 
to a graph with 4 vertices after applying a 1-step contraction ( all edges will remain black).

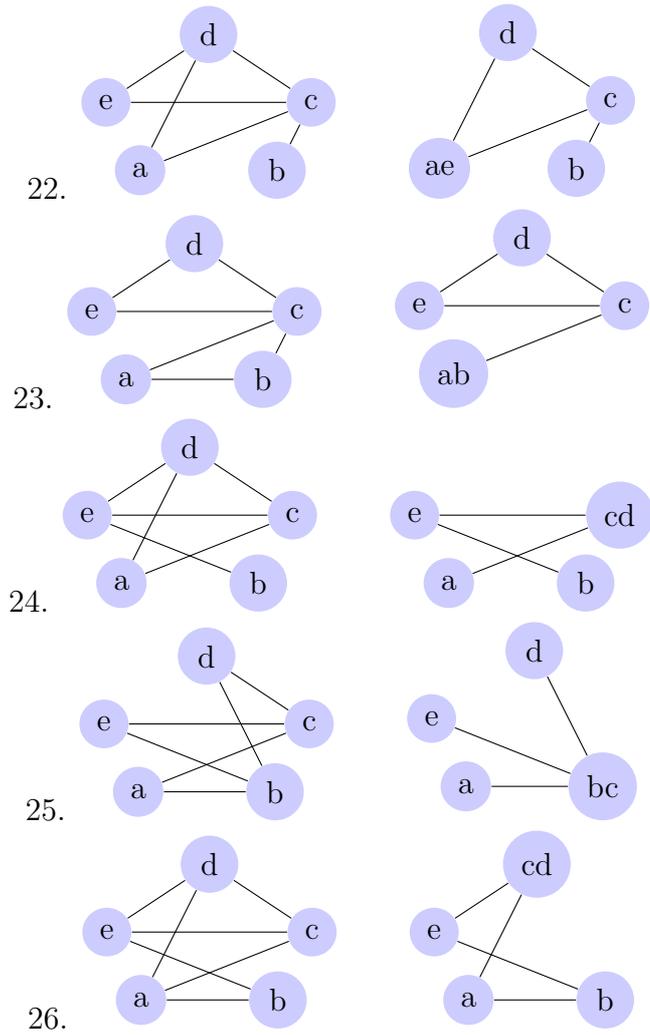
\begin{center}
22.
\begin{tikzpicture} 
[scale=.9,auto=center,every node/.style={circle,fill=blue!20}] 
\node (a1) at (-1,0) {a}; 
\node (a2) at (1,0) {b}; 
\node (a3) at (3/2,1) {c}; 
\node (a4) at (0,2) {d};
\node (a5) at (-3/2,1) {e};  

%\draw (a1) -- (a2);
\draw (a1) -- (a3);
\draw (a1) -- (a4);
%\draw (a1) -- (a5);
\draw (a2) -- (a3);
%\draw (a2) -- (a4);
%\draw (a2) -- (a5);
\draw (a3) -- (a4);
\draw (a3) -- (a5);
\draw (a4) -- (a5);

\end{tikzpicture}
\qquad
\begin{tikzpicture} 
[scale=.9,auto=center,every node/.style={circle,fill=blue!20}] 
\node (a1) at (-1,0) {ae}; 
\node (a2) at (1,0) {b}; 
\node (a3) at (3/2,1) {c}; 
\node (a4) at (0,2) {d};
%\node (a5) at (-3/2,1) {e};  

%\draw (a1) -- (a2);
\draw (a1) -- (a3);
\draw (a1) -- (a4);
%\draw (a1) -- (a5);
\draw (a2) -- (a3);
%\draw (a2) -- (a4);
%\draw (a2) -- (a5);
\draw (a3) -- (a4);
%\draw (a3) -- (a5);
%\draw (a4) -- (a5);
\end{tikzpicture}
\end{center}

\begin{center}
23.
\begin{tikzpicture} 
[scale=.9,auto=center,every node/.style={circle,fill=blue!20}] 
\node (a1) at (-1,0) {a}; 
\node (a2) at (1,0) {b}; 
\node (a3) at (3/2,1) {c}; 
\node (a4) at (0,2) {d};
\node (a5) at (-3/2,1) {e};  

\draw (a1) -- (a2);
\draw (a1) -- (a3);
%\draw (a1) -- (a4);
%\draw (a1) -- (a5);
\draw (a2) -- (a3);
%\draw (a2) -- (a4);
%\draw (a2) -- (a5);
\draw (a3) -- (a4);
\draw (a3) -- (a5);
\draw (a4) -- (a5);

\end{tikzpicture}
\qquad
\begin{tikzpicture} 
[scale=.9,auto=center,every node/.style={circle,fill=blue!20}] 
\node (a1) at (-1,0) {ab}; 
%\node (a2) at (1,0) {b}; 
\node (a3) at (3/2,1) {c}; 
\node (a4) at (0,2) {d};
\node (a5) at (-3/2,1) {e};  

%\draw (a1) -- (a2);
\draw (a1) -- (a3);
%\draw (a1) -- (a4);
%\draw (a1) -- (a5);
%\draw (a2) -- (a3);
%\draw (a2) -- (a4);
%\draw (a2) -- (a5);
\draw (a3) -- (a4);
\draw (a3) -- (a5);
\draw (a4) -- (a5);

\end{tikzpicture}
\end{center}

\begin{center}
24.
\begin{tikzpicture} 
[scale=.9,auto=center,every node/.style={circle,fill=blue!20}] 
\node (a1) at (-1,0) {a}; 
\node (a2) at (1,0) {b}; 
\node (a3) at (3/2,1) {c}; 
\node (a4) at (0,2) {d};
\node (a5) at (-3/2,1) {e};  

%\draw (a1) -- (a2);
\draw (a1) -- (a3);
\draw (a1) -- (a4);
%\draw (a1) -- (a5);
%\draw (a2) -- (a3);
%\draw (a2) -- (a4);
\draw (a2) -- (a5);
\draw (a3) -- (a4);
\draw (a3) -- (a5);
\draw (a4) -- (a5);

\end{tikzpicture}
\qquad
\begin{tikzpicture} 
[scale=.9,auto=center,every node/.style={circle,fill=blue!20}] 
\node (a1) at (-1,0) {a}; 
\node (a2) at (1,0) {b}; 
\node (a3) at (3/2,1) {cd}; 
%\node (a4) at (0,2) {d};
\node (a5) at (-3/2,1) {e};  

%\draw (a1) -- (a2);
\draw (a1) -- (a3);
%\draw (a1) -- (a4);
%\draw (a1) -- (a5);
%\draw (a2) -- (a3);
%\draw (a2) -- (a4);
\draw (a2) -- (a5);
%\draw (a3) -- (a4);
\draw (a3) -- (a5);
%\draw (a4) -- (a5);

\end{tikzpicture}
\end{center}

\begin{center}
25.
\begin{tikzpicture} 
[scale=.9,auto=center,every node/.style={circle,fill=blue!20}] 
\node (a1) at (-1,0) {a}; 
\node (a2) at (1,0) {b}; 
\node (a3) at (3/2,1) {c}; 
\node (a4) at (0,2) {d};
\node (a5) at (-3/2,1) {e};  

\draw (a1) -- (a2);
\draw (a1) -- (a3);
%\draw (a1) -- (a4);
%\draw (a1) -- (a5);
%\draw (a2) -- (a3);
\draw (a2) -- (a4);
\draw (a2) -- (a5);
\draw (a3) -- (a4);
\draw (a3) -- (a5);
%\draw (a4) -- (a5);

\end{tikzpicture}
\qquad
\begin{tikzpicture} 
[scale=.9,auto=center,every node/.style={circle,fill=blue!20}] 
\node (a1) at (-1,0) {a}; 
\node (a2) at (1,0) {bc}; 
%\node (a3) at (3/2,1) {c}; 
\node (a4) at (0,2) {d};
\node (a5) at (-3/2,1) {e};  

\draw (a1) -- (a2);
%\draw (a1) -- (a3);
%\draw (a1) -- (a4);
%\draw (a1) -- (a5);
%\draw (a2) -- (a3);
\draw (a2) -- (a4);
\draw (a2) -- (a5);
%\draw (a3) -- (a4);
%\draw (a3) -- (a5);
%\draw (a4) -- (a5);

\end{tikzpicture}
\end{center}

\begin{center}
26.
\begin{tikzpicture} 
[scale=.9,auto=center,every node/.style={circle,fill=blue!20}] 
\node (a1) at (-1,0) {a}; 
\node (a2) at (1,0) {b}; 
\node (a3) at (3/2,1) {c}; 
\node (a4) at (0,2) {d};
\node (a5) at (-3/2,1) {e};  

\draw (a1) -- (a2);
\draw (a1) -- (a3);
\draw (a1) -- (a4);
%\draw (a1) -- (a5);
%\draw (a2) -- (a3);
%\draw (a2) -- (a4);
\draw (a2) -- (a5);
\draw (a3) -- (a4);
\draw (a3) -- (a5);
\draw (a4) -- (a5);

\end{tikzpicture}
\qquad
\begin{tikzpicture} 
[scale=.9,auto=center,every node/.style={circle,fill=blue!20}] 
\node (a1) at (-1,0) {a}; 
\node (a2) at (1,0) {b}; 
%\node (a3) at (3/2,1) {c}; 
\node (a4) at (0,2) {cd};
\node (a5) at (-3/2,1) {e};  

\draw (a1) -- (a2);
%\draw (a1) -- (a3);
\draw (a1) -- (a4);
%\draw (a1) -- (a5);
%\draw (a2) -- (a3);
%\draw (a2) -- (a4);
\draw (a2) -- (a5);
%\draw (a3) -- (a4);
%\draw (a3) -- (a5);
\draw (a4) -- (a5);

\end{tikzpicture}
\end{center}
Since the graphs with 4 vertices have twin-width less than equal to 1, the graphs in the above-mentioned list will also have twin-width less than equal to 1. Next, we have a graph with twin-width less than equal to 2. 
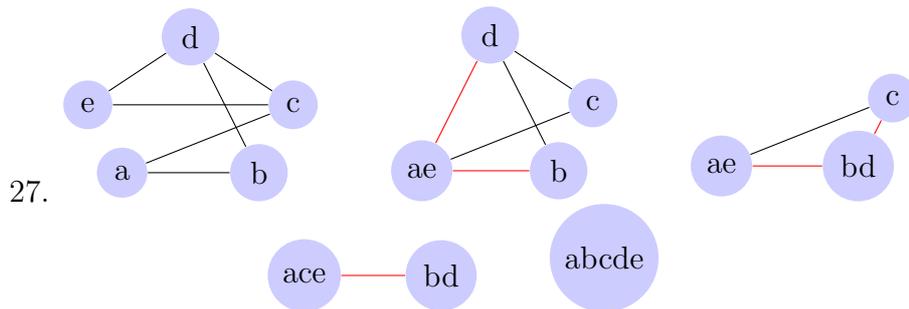
\begin{center}
27. 
\begin{tikzpicture} 
[scale=.9,auto=center,every node/.style={circle,fill=blue!20}] 
\node (a1) at (-1,0) {a}; 
\node (a2) at (1,0) {b}; 
\node (a3) at (3/2,1) {c}; 
\node (a4) at (0,2) {d};
\node (a5) at (-3/2,1) {e};  

\draw (a1) -- (a2);
\draw (a1) -- (a3);
%\draw (a1) -- (a4);
%\draw (a1) -- (a5);
%\draw (a2) -- (a3);
\draw (a2) -- (a4);
%\draw (a2) -- (a5);
\draw (a3) -- (a4);
\draw (a3) -- (a5);
\draw (a4) -- (a5);

\end{tikzpicture}
\qquad
\begin{tikzpicture} 
[scale=.9,auto=center,every node/.style={circle,fill=blue!20}] 
\node (a1) at (-1,0) {ae}; 
\node (a2) at (1,0) {b}; 
\node (a3) at (3/2,1) {c}; 
\node (a4) at (0,2) {d};
%\node (a5) at (-3/2,1) {e};  

\draw[red] (a1) -- (a2);
\draw (a1) -- (a3);
\draw[red] (a1) -- (a4);
%\draw (a1) -- (a5);
%\draw (a2) -- (a3);
\draw (a2) -- (a4);
%\draw (a2) -- (a5);
\draw (a3) -- (a4);
%\draw (a3) -- (a5);
%\draw (a4) -- (a5);

\end{tikzpicture}
\qquad
\begin{tikzpicture} 
[scale=.9,auto=center,every node/.style={circle,fill=blue!20}] 
\node (a1) at (-1,0) {ae}; 
\node (a2) at (1,0) {bd}; 
\node (a3) at (3/2,1) {c}; 
%\node (a4) at (0,2) {d};
%\node (a5) at (-3/2,1) {e};  

\draw[red] (a1) -- (a2);
\draw (a1) -- (a3);
\draw[red] (a2) -- (a3);
%\draw (a1) -- (a5);
%\draw (a2) -- (a3);
%\draw (a2) -- (a4);
%\draw (a2) -- (a5);
%\draw (a3) -- (a4);
%\draw (a3) -- (a5);
%\draw (a4) -- (a5);

\end{tikzpicture}
\qquad
\begin{tikzpicture} 
[scale=.9,auto=center,every node/.style={circle,fill=blue!20}] 
\node (a1) at (-1,0) {ace}; 
\node (a2) at (1,0) {bd}; 
%\node (a3) at (3/2,1) {c}; 
%\node (a4) at (0,2) {d};
%\node (a5) at (-3/2,1) {e};  

%\draw (a1) -- (a2);
\draw[red] (a1) -- (a2);
\end{tikzpicture}
\qquad
\begin{tikzpicture} 
[scale=.9,auto=center,every node/.style={circle,fill=blue!20}] 
\node (a1) at (0,0) {abcde}; 
\end{tikzpicture}
\end{center}

Now, we make a list of graphs whose complement graphs are disconnected. Since a graph has same twin-width as its complement graphs by Theorem \ref{complementgraph}, the graphs in this list will have twin-width less than equal to 1.

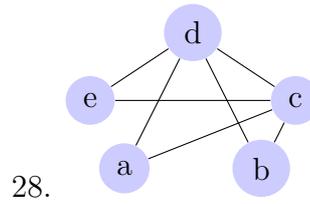
\begin{center}
28.
\begin{tikzpicture} 
[scale=.9,auto=center,every node/.style={circle,fill=blue!20}] 
\node (a1) at (-1,0) {a}; 
\node (a2) at (1,0) {b}; 
\node (a3) at (3/2,1) {c}; 
\node (a4) at (0,2) {d};
\node (a5) at (-3/2,1) {e};  

%\draw (a1) -- (a2);
\draw (a1) -- (a3);
\draw (a1) -- (a4);
%\draw (a1) -- (a5);
\draw (a2) -- (a3);
\draw (a2) -- (a4);
%\draw (a2) -- (a5);
\draw (a3) -- (a4);
\draw (a3) -- (a5);
\draw (a4) -- (a5);

\end{tikzpicture}
\end{center}

(The vertex `c' (or `d') is connected by edges with every other vertices. Therefore, the 
complement graph is disconnected.)

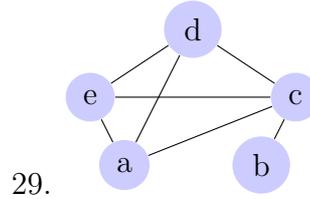
\begin{center}
29.
\begin{tikzpicture} 
[scale=.9,auto=center,every node/.style={circle,fill=blue!20}] 
\node (a1) at (-1,0) {a}; 
\node (a2) at (1,0) {b}; 
\node (a3) at (3/2,1) {c}; 
\node (a4) at (0,2) {d};
\node (a5) at (-3/2,1) {e};  

%\draw (a1) -- (a2);
\draw (a1) -- (a3);
\draw (a1) -- (a4);
\draw (a1) -- (a5);
\draw (a2) -- (a3);
%\draw (a2) -- (a4);
%\draw (a2) -- (a5);
\draw (a3) -- (a4);
\draw (a3) -- (a5);
\draw (a4) -- (a5);

\end{tikzpicture}
\end{center}
(The vertex `c' is connected by edges with every other vertices. Therefore, the complement graph is disconnected. )

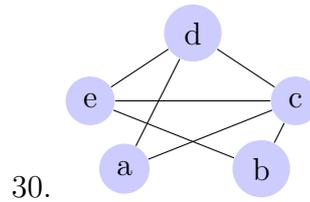
\begin{center}
30.
\begin{tikzpicture} 
[scale=.9,auto=center,every node/.style={circle,fill=blue!20}] 
\node (a1) at (-1,0) {a}; 
\node (a2) at (1,0) {b}; 
\node (a3) at (3/2,1) {c}; 
\node (a4) at (0,2) {d};
\node (a5) at (-3/2,1) {e};  

%\draw (a1) -- (a2);
\draw (a1) -- (a3);
\draw (a1) -- (a4);
%\draw (a1) -- (a5);
\draw (a2) -- (a3);
%\draw (a2) -- (a4);
\draw (a2) -- (a5);
\draw (a3) -- (a4);
\draw (a3) -- (a5);
\draw (a4) -- (a5);

\end{tikzpicture}
\end{center}
(The vertex `c' is connected by edges with every other vertices. Therefore, the complemented graph is 
disconnected. )

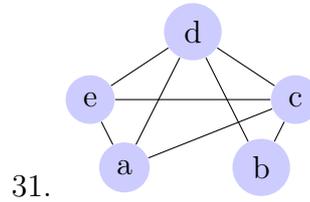
\begin{center}
31.
\begin{tikzpicture} 
[scale=.9,auto=center,every node/.style={circle,fill=blue!20}] 
\node (a1) at (-1,0) {a}; 
\node (a2) at (1,0) {b}; 
\node (a3) at (3/2,1) {c}; 
\node (a4) at (0,2) {d};
\node (a5) at (-3/2,1) {e};  

%\draw (a1) -- (a2);
\draw (a1) -- (a3);
\draw (a1) -- (a4);
\draw (a1) -- (a5);
\draw (a2) -- (a3);
\draw (a2) -- (a4);
%\draw (a2) -- (a5);
\draw (a3) -- (a4);
\draw (a3) -- (a5);
\draw (a4) -- (a5);

\end{tikzpicture}
\end{center}
(The vertex `c' (or `d') is connected by edges with other vertices. Therefore, the complement graph is disconnected. )

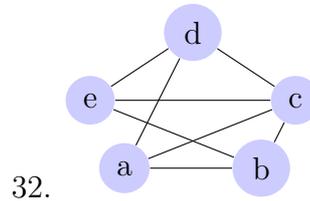
\begin{center}
32.
\begin{tikzpicture} 
[scale=.9,auto=center,every node/.style={circle,fill=blue!20}] 
\node (a1) at (-1,0) {a}; 
\node (a2) at (1,0) {b}; 
\node (a3) at (3/2,1) {c}; 
\node (a4) at (0,2) {d};
\node (a5) at (-3/2,1) {e};  

\draw (a1) -- (a2);
\draw (a1) -- (a3);
\draw (a1) -- (a4);
%\draw (a1) -- (a5);
\draw (a2) -- (a3);
%\draw (a2) -- (a4);
\draw (a2) -- (a5);
\draw (a3) -- (a4);
\draw (a3) -- (a5);
\draw (a4) -- (a5);

\end{tikzpicture}
\end{center}
(The vertex `c' is connected by edges with every other vertices. Therefore, the complement graph is disconnected.)

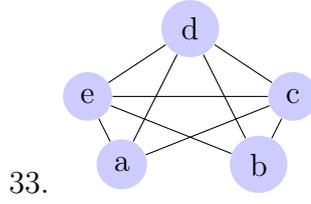
\begin{center}
33.
\begin{tikzpicture} 
[scale=.9,auto=center,every node/.style={circle,fill=blue!20}] 
\node (a1) at (-1,0) {a}; 
\node (a2) at (1,0) {b}; 
\node (a3) at (3/2,1) {c}; 
\node (a4) at (0,2) {d};
\node (a5) at (-3/2,1) {e};  

%\draw (a1) -- (a2);
\draw (a1) -- (a3);
\draw (a1) -- (a4);
\draw (a1) -- (a5);
\draw (a2) -- (a3);
\draw (a2) -- (a4);
\draw (a2) -- (a5);
\draw (a3) -- (a4);
\draw (a3) -- (a5);
\draw (a4) -- (a5);

\end{tikzpicture}
\end{center}

(The vertex `c' (or `d' or `e') is connected by edges with every other vertices. Therefore, the complement graph is disconnected. )

Finally, we are left with the following complete graph whose twin-width is zero. 

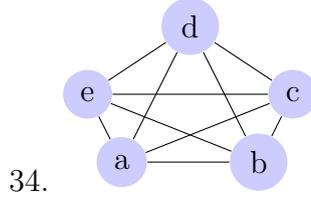
\begin{center}
34.
\begin{tikzpicture} 
[scale=.9,auto=center,every node/.style={circle,fill=blue!20}] 
\node (a1) at (-1,0) {a}; 
\node (a2) at (1,0) {b}; 
\node (a3) at (3/2,1) {c}; 
\node (a4) at (0,2) {d};
\node (a5) at (-3/2,1) {e};  

\draw (a1) -- (a2);
\draw (a1) -- (a3);
\draw (a1) -- (a4);
\draw (a1) -- (a5);
\draw (a2) -- (a3);
\draw (a2) -- (a4);
\draw (a2) -- (a5);
\draw (a3) -- (a4);
\draw (a3) -- (a5);
\draw (a4) -- (a5);

\end{tikzpicture}
\end{center}

Hence, we have our theorem
\hfill\(\Box\)

\section{Upper bound of twin-width of King's graph and Rook's graph}

In this section, we prove Theorem \ref{King'sgraph} and Theorem \ref{Rook'sgraph}. Before going
into the proof, we introduce two types of graph products Strong Product and Cartesian Product which will be required in the proofs of the theorems. The \textit{strong product} of graphs $G$ and $H$, denoted by $G\boxtimes H$, is a graph such that the vertex set of $G\boxtimes H$ is the Cartesian product $V(G) \times V(H)$, where $V(G)$ and $V(H)$ are the set of vertices of $G$ and $H$, respectively; and distinct vertices $(u,u')$ and $(v,v')$ are adjacent in $G\boxtimes H$ if and only if $u = v$ and $u'$ is adjacent to $v'$, or $u' = v'$ and $u$ is adjacent to $v$, or
$u$ is adjacent to $v$ and $u'$ is adjacent to $v'$. On the other hand, the \textit{Cartesian product} of graphs $G$ and $H$, denoted by $G\square H$, is a graph such that the vertex set of $G\square H$ is the Cartesian product $V(G) \times V(H)$; and two vertices $(u,u')$ and $(v,v')$ are adjacent in $G\square H$ if and only if either $u = v$ and $u'$ is adjacent to $v'$ in $H$, or $u'$ = $v'$ and $u$ is adjacent to $v$ in G. However, we now prove Theorem \ref{King'sgraph}. 

\vspace{5mm}

\textbf{Proof of Theorem \ref{King'sgraph}:} The $(n\times m)$-King's graph, is a strong product of two path graphs $P_n$ and $P_m$.  
%We draw the $(8\times 8)$-King's graph below.  
%\begin{figure}[!ht]
%\begin{center}
%\includegraphics[width=0.3\textwidth,natwidth=610,natheight=642]{Kinggraph.jpg}
%\end{center}
%\end{figure}
We obtain from \cite{BGKTW21} (Theorem 9) that
$$tww(G\boxtimes H)\leq max \{tww(G)(\Delta(H) + 1) +2 \Delta(H), tww(H) + \Delta(H)\} .$$ 

Since the twin-width of a path graph is less than equal to $1$ and $\Delta(P_m)=2$, we obtain that the twin-width of a King's graph with $n$-vertices is less than equal to $7$.
\hfill\(\Box\)

Now, we prove Theorem \ref{Rook'sgraph}. 

\vspace{5mm}

\textbf{Proof of Theorem \ref{Rook'sgraph}:}  The $(n\times m)$-Rook's graph is a Cartesian product of two Complete graphs $K_n$ and $K_m$.
%We draw the $(8\times 8)$-Rook's graph  below. 
%\begin{figure}[!ht]
%\begin{center}
%\includegraphics[width=0.3\textwidth,natwidth=610,natheight=642]{Rookgraph.jpg}
%\end{center}
%\end{figure}
We obtain from \cite{PS22} (Theorem 3.1) that for any graphs $G$ and $H$, 
$$tww(G\square H)\leq max\{ tww(G) + \Delta(H), tww(H)\} + \Delta(H).$$
Since the twin-width of a complete graph is zero and $\Delta(K_m)$ is $(m-1)$, we have the 
twin-width of $(n\times m)$-Rook's graph is less than equal to  $2(m-1)$. 
\hfill\(\Box\)

\end{document}